\documentclass[11pt]{amsart}
\usepackage[utf8]{inputenc}

\usepackage{subfiles}

\usepackage{amsmath}
\usepackage{amssymb}
\usepackage{amsthm}
\usepackage{thmtools}
\usepackage{thm-restate}
\usepackage[margin=1in]{geometry}
\usepackage{cleveref}

\numberwithin{equation}{section}

\DeclareMathOperator{\ord}{ord}
\DeclareMathOperator{\SL}{SL}

\theoremstyle{definition}
\newtheorem{theorem}{Theorem}[section]
\newtheorem*{theorem*}{Theorem}

\newtheorem*{example*}{Example}
\newtheorem{lemma}[theorem]{Lemma}
\newtheorem*{lemma*}{Lemma}

\newtheorem*{corollary*}{Corollary}

\newtheorem*{definition*}{Definition}
\newtheorem{proposition}[theorem]{Proposition}
\newtheorem*{proposition*}{Proposition}

\newtheorem*{remark*}{Remark}

\newcommand{\mD}{\mathcal{D}}
\newcommand{\mO}{\mathcal{O}}

\newcommand{\Z}{\mathbb{Z}}
\newcommand{\Q}{\mathbb{Q}}

\newcommand{\C}{\mathbb{C}}

\newcommand{\tm}{\times}

\newcommand{\bmat}[1]{\begin{bmatrix} #1 \end{bmatrix}}

\renewcommand{\aa}{\alpha}
\newcommand{\bb}{\beta}
\newcommand{\cc}{\gamma}
\newcommand{\dd}{\delta}
\newcommand{\om}{\omega}
\newcommand{\Th}{\Theta}

\newcommand{\ps}[1]{\left(#1\right)}
\newcommand{\floor}[1]{\left\lfloor #1 \right\rfloor}
\newcommand{\abs}[1]{\left| #1 \right|}
\newcommand{\leg}[2]{\left(\frac{#1}{#2}\right)}

\newcommand{\poc}[1]{(#1; #1)_{\infty}}

\title{Ramanujan Congruences for Fractional Partition Functions}
\author{Erin Bevilacqua, Kapil Chandran, and Yunseo Choi}

\begin{document}
\maketitle

\begin{abstract}
For rational $\alpha$, the fractional partition functions $p_\alpha(n)$ are given by the coefficients of the generating function $(q;q)^\alpha_\infty$.
When $\alpha=-1$, one obtains the usual partition function. Congruences of the form $p(\ell n + c)\equiv 0 \pmod{\ell}$ for a prime $\ell$ and integer $c$ were studied by Ramanujan. Such congruences exist only for $\ell\in\{5,7,11\}.$
Chan and Wang \cite{chan} recently studied congruences for the fractional partition functions and gave several infinite families of congruences using identities of the Dedekind eta-function. Following their work, we use the theory of non-ordinary primes to find a general framework that characterizes congruences modulo any integer. 
This allows us to prove new congruences such as $p_\frac{57}{61}(17^2n-3)\equiv 0 \pmod{17^2}$. 
\end{abstract}

\section{Introduction}
\label{intro}


A \textit{partition} of a nonnegative integer $n$ is a non-increasing sequence of positive integers that sum to $n$.
For every nonnegative integer $n$, denote by $p(n)$ the number of partitions of $n$. 
The generating function for $p(n)$ was shown by Euler to be
\[ P(q):= \sum_{n=0}^{\infty} p(n)q^n = \frac{1}{\poc{q}}, \]
where $\poc{q} := \prod_{n=1}^{\infty} (1-q^{n})$ is the $q$-Pochhammer symbol.
Ramanujan observed the following congruences for every nonnegative integer $n$:
\begin{align*}
    p(5n+4) &\equiv 0 \pmod{5}, \\
    p(7n+5) &\equiv 0 \pmod{7}, \\
    p(11n+6) &\equiv 0 \pmod{11}.
\end{align*}
Ramanujan proved the first two congruences, while Atkin proved the third \cite{web}. Ahlgren and Boylan \cite{ahlgren} later proved that congruences of this form, in which the modulus of the congruence and common difference of the arithmetic progression are equal to the same prime, are limited to these three identified by Ramanujan.

Beyond these congruences, Ramanujan conjectured that similar congruences hold for all powers of $5$, $7$, and $11.$
An infinite class of congruences modulo powers of $5$, $7$, and $11$ has since been proven \cite{web}, generalizing the previous congruences and confirming Ramanujan's intuition. These higher congruences can be described as follows.
For a positive integer $k$ and a prime $\ell$, let $\dd_{\ell, k}$ be the unique integer in the set $\{0, 1, 2, \ldots, \ell^k - 1\}$ such that $\delta_{\ell,k} \equiv 1/24 \pmod{\ell^k}.$ Then for every nonnegative integer $n$,
\begin{align*}
    p(5^kn + \delta_{5,k}) &\equiv 0 \pmod{5^k},\\
    p(7^kn + \delta_{7,k}) &\equiv 0 \pmod{7^{\floor{k/2}+1}},\\
    p(11^kn + \delta_{11,k}) &\equiv 0 \pmod{11^k}.
\end{align*}

These are the only congruences of this shape for the partition function, but many more congruences exist when we relax the condition that the modulus of the congruence and the common difference of the arithmetic progression must both be powers of the same prime $\ell$. Generalizing in this way, Ono \cite{annals} showed how to find congruences for the partition function modulo any prime $\ell \ge 5$. The work of Ahlgren and Ono \cite{ahlono} expands upon this result, guaranteeing that for every integer $L$ coprime to $6$, there exist integers $A$ and $B$ such that $p(An+B) \equiv 0 \pmod{L}$ for all nonnegative integers $n$. 
An example of such a congruence is \[ p(4063467631n+30064597)\equiv 0 \pmod{31}. \]

As an extension of the theory of congruences for the partition function, Chan and Wang \cite{chan} recently explored congruences for the fractional partition functions.
For rational $\aa$, which we write as $\aa = a/b$ for coprime integers $a, b$ with $b \ge 1$ throughout this paper, the fractional partition function $p_{\alpha}\colon \Z \to \Q$ is defined for $n \ge 0$ by the power series
\begin{equation}\label{fracp}
   P_{\alpha}(q):= {\poc{q}^{\alpha}}
   =: \sum_{n=0}^{\infty} p_{\alpha}(n)q^n.
\end{equation}
We set $p_\aa(n) := 0$ for $n < 0.$ 
Chan and Wang \cite{chan} showed (see \Cref{definitions}) that the numbers $p_\aa(n)$ are $\ell$-integral for any prime $\ell \nmid b$. 
For this reason, it is fruitful to study congruences for $p_\aa$ modulo $\ell^r$ for any prime $\ell \nmid b.$ We study \textit{$\ell^r$-balanced} congruences, i.e. those of the form,
\begin{equation}
p_\aa(\ell^r n + c) \equiv 0 \pmod{\ell^r}
\end{equation}
for all $n$, where $c, r$ are integers with $r \ge 1.$
In other words, an $\ell^r$-balanced congruence is one in which the modulus of the congruence and the common difference of the arithmetic progression are both equal to the prime power $\ell^r$. 

Balanced congruences have been studied for integral powers of the partition function. In \cite{K+O}, Kiming and Olsson identified necessary conditions for $p_\aa$ to admit an $\ell$-balanced congruence for integer $\aa.$ Their main result (Theorem 1 of \cite{K+O}) can be stated as follows. Suppose that $\aa$ is an integer and $\ell \ge 5$ is a prime such that $\aa \not\equiv 1, 3 \pmod{\ell}.$ If $p_\aa(\ell n + c) \equiv 0 \pmod{\ell}$  for all $n$, then $24c + \aa \equiv 0 \pmod{\ell}.$ This result explains the significance of the residue class $1/24 \pmod{\ell}$ in Ramanujan's congruences for the partition function.
Our next theorem is an analog of Kiming and Olsson's result for fractional $\alpha.$ It shows that when $\alpha \not \equiv 1, 3 \pmod{\ell}$, there is generically at most one residue class modulo $\ell$ corresponding to an $\ell$-balanced congruence for $p_\aa.$

\begin{restatable}{theorem}{kiming}
\label{kiming}
Let $\aa = a/b$, and let $\ell \ge 5$ be a prime not dividing $b$ such that $\aa \not\equiv 1, 3 \pmod{\ell}.$ If \[ p_\aa(\ell n + c) \equiv 0 \pmod{\ell} \] for all $n,$ then $24c + \aa \equiv 0 \pmod{\ell}.$
\end{restatable}

The cases where $\alpha\equiv 1,3\pmod{\ell}$ yield far more $\ell$-balanced congruences. These two cases are distinguished by two famous lacunary Fourier expansion identities for powers of the Dedekind eta-function ${\eta(z) := q^{1/24} \poc{q}}$, where $q := e^{2\pi iz}.$ Recall that a power series $\sum_{n = 0}^{\infty} a(n)q^n$ is called \textit{lacunary} if 
\[ \lim_{N \to \infty} \frac{\#\{n \le N : a(n) = 0\}}{N} = 1. \] 
In \cite{lacunary}, Serre studied and classified which even positive powers of the eta-function are lacunary.
Explicitly, he found that for even integer $r > 0$, $\eta^r$ is lacunary if and only if $r \in \{2, 4, 6, 8, 10, 14, 26\}.$
In the case when $r$ is odd, we have the following classical formulae due to Euler and Jacobi, respectively:
\begin{align}
\label{euler}
    \eta(24z) &= \sum_{n = 1}^{\infty} \leg{12}{n}q^{n^2}, \\
\label{jacobi}
    \eta(8z)^3 &= \sum_{n = 1}^{\infty} \leg{-4}{n}nq^{n^2},
\end{align}
where $\leg{\cdot}{\cdot}$ is the Kronecker symbol. These formulae show that $\eta(24z)$ and $\eta(8z)^3$ are lacunary. Chan and Wang \cite{chan} used this lacunarity to obtain $\ell$-balanced congruences for $p_\aa$ when $\aa \equiv 1, 3 \pmod{\ell}.$ Our second result is an extension of Chan and Wang's work to include $\ell^r$-balanced congruences for $r > 1.$

\begin{restatable}{theorem}{thmtwo}
\label{thm2}
Let $\aa = a/b$, let $\ell$ be a prime not dividing $b$, and let $c$ be an integer.
\begin{itemize}
    \item[(i)] If $\leg{24c + 1}{\ell} = -1$, then for all $n$ we have \[ p_{\aa}(\ell n + c) \equiv 0 \pmod{\ell^{\ord_\ell(\aa - 1)}}. \]
    
    \item[(ii)] If $\leg{8c + 1}{\ell} \ne 1$, then for all $n$ we have \[ p_{\aa}(\ell n + c) \equiv 0 \pmod{\ell^{\ord_\ell(\aa - 3)}}. \] 
\end{itemize}
\end{restatable}

\begin{remark*}
\label{thm3}
By the Chinese remainder theorem, a simple consequence of \Cref{thm2} is that for any positive integer $L$, there exists $\aa \in \Q$ such that for all $n$,
\[ p_{\alpha}(Ln + c) \equiv 0 \pmod{L}. \]
\end{remark*}

In addition to the exceptional closed formulas for $\eta$ and $\eta^3$, there are certain even powers of the eta-function that are well known to share the property of lacunarity arising from the theory of complex multiplication.
By exploiting these even powers mentioned above, Chan and Wang \cite{chan} found $\ell$-balanced congruences for $p_\alpha$ when $\alpha \equiv 4,6,8,10,14,26 \pmod{\ell}$.
Their methods produce a large number of $\ell$-balanced congruences, expanding the results of previous literature which only presented a handful. 

We seek to obtain congruences through more general properties of modular forms that are not constrained to the finite set of even powers of the eta-function that are lacunary. 
To do so, we make use of a different phenomenon which is universal for modular forms, thereby establishing a general framework for describing all known $\ell$-balanced congruences for $p_\aa.$ 
The advantage to this method is that we can obtain $\ell$-balanced congruences for $p_\aa$ when $\aa$ is not restricted to a finite set of residue classes modulo $\ell.$
Furthermore, our theory enables us to extend our results to obtain $\ell^r$-balanced congruences for $r > 1.$

In order to exploit this phenomenon, we need the following definition. We say that a prime \textit{$\ell$ is good for $\alpha = a/b$ with parameter $k$} if $\ell \nmid b$ and $k$ is a positive integer such that
\begin{flalign}
\label{cond1}
    & \ell \mid (24k - \aa),&\\
\label{cond2}
    & (\ell - 1) \mid (12k - m) \text{ for some } m\in\{4,6,8,10,14\}, \text{ and}&\\
\label{cond3}
    & \ell\nmid N_{12k}(\mD_{12k}) \text{ where } \mD_{12k} \text{ is the Hecke determinant for } S_{12k} \text{ (see \Cref{lnonordinary}). }&
\end{flalign}
We say that \textit{$\ell$ is good for $\aa$} if there exists some such parameter $k.$

\begin{restatable}{theorem}{nonordthm}
\label{nonordthm}
If $\ell$ is good for $\aa$ with parameter $k$ and $v \leq \ord_\ell(24k - \aa)$ is a positive integer, then for all $n$, we have \[ p_{\aa}(\ell^v n - k) \equiv 0 \pmod{\ell^v}. \]
\end{restatable}

\begin{remark*}
If a prime $\ell$ is good for $\alpha$ with parameter $k$, then $k < \ell$ (see \Cref{kltell}) and $\ell \le \max\{\abs{26b - a}, \abs{14b + a}\}$ for all but finitely many $\aa$ (see \Cref{finitegood}). Thus, we conclude that for generic $\aa$, there are only finitely many primes $\ell$ good for $\alpha$.
We note that the $\alpha$ for which we cannot derive a bound on $\ell$ as in \Cref{finitegood} is $\alpha \in \{0, 2, 4, 6, 8, 10, 14, 26\}$, aligning with the lacunary even powers of the eta-function classified by Serre. 
\end{remark*}

\Cref{nonordthm} is a consequence of the $\ell$-non-ordinarity of Hecke eigenforms (see \Cref{lnonordinary}). For a prime $\ell$, we say that a modular form $f(z) = \sum_{n = 0}^{\infty} a(n)q^n \in M_k \cap \mO_L[[q]]$, where $\mO_L$ is the ring of integers of a number field $L$, is \textit{$\ell$-non-ordinary} if $a(\ell) \equiv 0 \pmod{\ell\mO_L}.$ The notion of $\ell$-non-ordinarity is a slight strengthing of the notion of non-ordinarity at $\ell$ studied by Jin, Ma, Ono \cite{non-ord}, and others. If $f$ is a normalized Hecke eigenform, $\ell$-non-ordinarity is equivalent to $f$ being annihilated by the Hecke operator $T_\ell$ modulo $\ell$, that is, $f\mid T_\ell \equiv 0 \pmod{\ell\mO_L}.$
The second condition \eqref{cond2} for $\ell$ being good for $\aa$ is a technical criterion that guarantees the $\ell$-non-ordinarity of every normalized eigenform in $S_{12k}$ via a result of Jin, Ma, and Ono \cite{non-ord} (see \Cref{lnonordinary}). The third condition \eqref{cond3} ensures that the $\ell$-non-ordinarity of the normalized eigenforms extends to $\ell$-non-ordinarity for suitable linear combinations of these eigenforms. The first condition \eqref{cond1} ensures that the conditions of \Cref{kiming} are met. 

\Cref{nonordthm} represents a complete classification of the situations for which congruences of the fractional partition functions arise ``automatically'' because every normalized eigenform in $S_{12k}$ is $\ell$-non-ordinary. 
Combining results from \Cref{thm2} and \Cref{nonordthm}, we can deduce all congruences achieved by the main theorem of Chan and Wang (Theorem 1.2 of \cite{chan}) up to the third condition of $\ell$ being good for $\aa$ with parameter $k.$ 

\begin{example*}
Consider Theorem 1.2 (5) of \cite{chan}, which equivalently states that if $\ell \equiv 11 \pmod{12}$ and $\aa \equiv 26 \pmod{\ell}$, then for any integer $c$ satisfying $24c + 26 \equiv 0 \pmod{\ell}$, we have
\begin{equation}
\label{chanexample}
    p_\aa(\ell n + c) \equiv 0 \pmod{\ell}
\end{equation} 
for all $n.$ We can prove this result using \Cref{nonordthm} up to condition \eqref{cond3}. To do so, we choose $1 \le k < \ell$ such that $k \equiv -c \pmod{\ell}.$ We verify the first two conditions of $\ell$ being good for $\aa$ with parameter $k.$ Clearly $24k - \aa \equiv -24c - 26 \equiv 0 \pmod{\ell}$, which is \eqref{cond1}. Because $\ell \equiv 11 \pmod{12}$, it can be checked that $12k - 13 = \ell.$ Therefore $(\ell - 1) \mid (12k - 14)$, which is \eqref{cond2}. If \eqref{cond3} also holds, we can conclude by \Cref{nonordthm} that $p_\aa(\ell n - k) \equiv 0 \pmod{\ell}$ for all $n$, which is equivalent to \eqref{chanexample}.
\end{example*}

\begin{example*}
We use \Cref{nonordthm} to recover the $\ell$-balanced congruences for the usual partition function.
By \Cref{kltell} and \Cref{finitegood}, if a prime $\ell$ is good for $\aa = -1$ with parameter $k$, then $k < \ell$ and $\ell \le 27.$
One can check that the only values of $\ell \le 27$ which are good for $\aa = -1$ are $\ell \in \{5, 7, 11\}$, aligning with the results of Ahlgren and Boylan \cite{ahlgren}. For example, $\ell = 11$ is good for $\aa = -1$ with parameter $k = 5$ because $11 \nmid N_{60}(\mD_{60}).$  To compute the norm of the Hecke determinant $N_{60}(\mD_{60})$, we use the function \texttt{CuspForms(1,60).hecke\_algebra().discriminant()} on CoCalc. 
\end{example*}

\begin{example*}
Let $\ell=367$ and $\alpha=3/37$. Then, $\ell$ is good for $\alpha$ with parameter $k = 31$ because $367 \nmid N_{372}(\mD_{372}).$ We conclude from \Cref{nonordthm} that for all $n$, \[ p_\frac{3}{37}(367n-31) \equiv 0 \pmod{367}. \]
\end{example*}



\begin{example*}
We show that the set of $v$ identified by \Cref{nonordthm} is sharp for certain choices of $\alpha$ and $\ell$.
Let $\ell = 17$ and $\aa = 57/61$.
Then, $\ell$ is good for $\alpha$ with parameter $k = 3$ because $17 \nmid N_{36}(\mD_{36}).$
Since $\ord_{17}(24 \cdot 3 - \frac{57}{61}) = 2$, we conclude from \Cref{nonordthm} that for all $n$,
\begin{align*}
    p_\frac{57}{61}(17 n - 3) &\equiv 0 \pmod{17}, \\
    p_\frac{57}{61}(17^2 n - 3) &\equiv 0 \pmod{17^2}.
\end{align*}
However, we do not have an analogous $17^3$-balanced congruence because \[ p_\frac{57}{61}(17^3 - 3) = p_\frac{57}{61}(4910) \equiv \frac{1445}{2052} \not\equiv 0 \pmod{17^3}. \]
\end{example*}

\begin{remark*}
As mentioned above, \Cref{nonordthm} follows from from the $\ell$-non-ordinarity of Hecke eigenforms, with condition \eqref{cond2} ensuring the $\ell$-non-ordinarity of all normalized eigenforms in $S_{12k}$.
However, the prime $\ell=2$ is rarely good for $\alpha$ with parameter $k$, because condition \eqref{cond3} $2 \nmid N_{12k}(\mD_{12k})$ does not hold for any $k > 1$ (see \Cref{kltell}). While we do not do so here, it is possible to obtain a larger class of congruences for $p_\aa$ modulo powers of $2$ using Serre's theory of Hecke nilpotency \cite{web}.
By repeatedly applying the Hecke operator $T_2$, one can obtain congruences for the Fourier coefficients of a given power of $\Delta$ modulo a power of $2$ following similar steps as in the proof of \Cref{nonordthm}.


\end{remark*}

We now return to our study of necessary conditions for fractional partition function congruences. For a given $\aa$, we want to restrict the set of primes $\ell$ for which $p_\aa$ can admit an $\ell$-balanced congruence. 
Ahlgren and Boylan \cite{ahlgren} showed that for the usual partition function, $\ell$-balanced congruences exist only for $\ell \in \{5,7,11\}.$
Boylan \cite{boylan} later extended this result to show that for a negative odd integer $\aa$, there are finitely many primes $\ell$ for which $p_\aa$ admits an $\ell$-balanced congruence.
Using the same framework as Boylan, we extend this result to even $\aa < 0$ and odd $\aa > 3.$ 

\begin{restatable}{theorem}{finitel}
\label{finitel}
Let $\aa$ be an integer that is either even and $< 0$ or odd and $> 3.$ If \[ p_\aa(\ell n - \dd_\ell) \equiv 0 \pmod{\ell} \] for all $n,$ then $\ell \le |\aa| + 4.$ In particular, $p_\aa$ admits finitely many $\ell$-balanced congruences.
\end{restatable}


\begin{remark*}
It is unsurprising that this theorem does not extend to even positive $\alpha$.
When $\aa = 24$, for example, congruences for $p_\aa$ are equivalent to congruences for the Delta-function. According to \Cref{extendpowers}, congruences for $\Delta$ would be implied by $\ell$-non-ordinarity of $\Delta$ at different primes $\ell.$
Whether $\Delta$ is non-ordinary at infinitely many primes is a famous open problem \cite{story} that has been studied extensively.
\end{remark*}

The methods used by Ahlgren and Boylan \cite{ahlgren}, and later by Boylan \cite{boylan}, do not extend in the same way to all rational $\alpha$.
Instead, we obtain a a residue class restriction on the primes $\ell$ for which $p_\aa$ can admit an $\ell$-balanced congruence. 
To state our result, it is convenient to introduce the following notation. If $m$ is a positive integer and $\bb$ is a rational number whose denominator in lowest terms is coprime to $m$, we denote by $\Psi_m(\bb)$ the unique integer in the set $\{0, 1, \ldots, m - 1\}$ that is congruent to $\bb$ modulo $m.$

\begin{restatable}{theorem}{fiiltrationcor}
\label{filtrationcor}
Let $\aa = a/b$ be a rational number that is not an even integer $\geq 0$. If $\ell \ge |a| + 5b$ is a prime for which $p_\aa$ admits an $\ell$-balanced congruence, then \[ \Psi_{2b}\ps{\frac{a}{\Psi_{2b}(\ell)}} \ge b. \]
\end{restatable}



\begin{remark*}
This theorem shows that for sufficiently large $\ell$, half of the primes classified by residue class modulo $2b$ cannot be the modulus of a balanced congruence for $p_\aa$ due to Dirichlet's Theorem. 
\end{remark*}

The remainder of this paper is structured as follows. In \Cref{prelim}, we give an overview of properties of modular forms relevant to our research. We also provide some basic results from other papers and prove several lemmas required for our main theorems. In Section 3 we prove the main results. 

\section{Preliminaries}
\label{prelim}

In \Cref{definitions}, we discuss elementary definitions and facts about modular forms and introduce some basic results from Chan and Wang \cite{chan} that motivate our work.
In \Cref{lnonordinary}, we introduce several lemmas about $\ell$-non-ordinary primes.
Finally, in \Cref{modulo}, we introduce some results on Serre filtrations and modular forms modulo $\ell$ for prime $\ell.$

\subsection{Modular forms and preliminary results}
\label{definitions} 
The following facts about modular forms are well-known and can be found in any standard text, such as \cite{web}.
For an integer $k$, denote by $M_k$ (resp. $S_k$) the $\C$-vector space of holomorphic modular forms (resp. cusp forms) of weight $k$ on $\SL_2(\Z).$

Recall that for any positive integer $m$, the $m$th Hecke operator $T_{m, k}$ (often abbreviated to $T_m$ when the weight is clear) is an endomorphism on $M_k.$ Its action on a Fourier expansion $f(z) = \sum_{n = 0}^{\infty} a(n)q^n$ is given by the formula 
\begin{equation*}
    f(z) \mid T_{m,k} = \sum_{n=0}^{\infty} \ps{\sum_{d\mid (m,n)}
    d^{k-1}a(mn/d^2)}q^n.
\end{equation*}
For $m = \ell$ prime, this reduces to \[ f(z)\mid T_{\ell, k} := \sum_{n=0}^{\infty} \ps{a(\ell n) + \ell^{k-1}a(n/\ell)}q^n. \] We call a modular form $f(z) \in S_k$ a \textit{Hecke eigenform} if it is an eigenvector of $T_{m, k}$ for all $m \ge 1$, i.e. there exist $\lambda(m) \in \C$ such that \[ f(z)\mid T_{m,k} = \lambda(m)f(z). \] We say that a Hecke eigenform $f = \sum_{n = 0}^{\infty} a(n)q^n$ is \textit{normalized} if $a(1) = 1.$ 

We make extensive use of modular forms throughout this paper to prove congruences for the fractional partition functions. As a starting point, we require two key results from Chan and Wang's work \cite{chan}. The first result (Theorem 1.1 of \cite{chan}) identifies which congruences are meaningful to study.
\begin{theorem}\label{lintegral}
Let $\alpha = a/b.$ Then the denominator of $p_\aa(n)$ when written in lowest terms is given by
\[ \text{denom}(p_{\alpha}(n))=b^n\prod_{p\mid b}p^{\ord_p(n!)}. \]
\end{theorem}
From this theorem, we conclude that for a given rational number $\alpha$, we can study congruences for $p_\aa$ modulo $L$ whenever $\gcd(L, b) = 1.$ The second result that we require is a technical lemma (Lemma 2.1 of \cite{chan}) that is a consequence of applying the Frobenius endomorphism. This lemma allows us to move exponents through $q$-Pochhammer symbols, which is a crucial step in several proofs.
\begin{lemma}\label{movel}
Let $\aa = a/b.$ Let $\ell$ be a prime not dividing $b.$ Then for any $r \ge 1$, 
\begin{align*}
    \poc{q}^{\ell^r\aa} \equiv \poc{q^\ell}^{\ell^{r - 1}\aa} \pmod{\ell^r}.
\end{align*}
\end{lemma}

\subsection{$\ell$-non-ordinary primes}
\label{lnonordinary}
We now introduce the related notions of non-ordinarity at $\ell$ and $\ell$-non-ordinarity. Throughout this subsection, we denote by $L$ a number field with ring of integers $\mO_L.$
Recall that a modular form ${f(z)=\sum_{n=0}^\infty a(n)q^n \in M_k \cap \mO_L[[q]]}$ is said to be \textit{non-ordinary} at $\ell$ if there exists a prime ideal $\mathfrak{l}$ above $\ell$ such that \[ a(\ell)\equiv 0 \pmod{\mathfrak{l}}. \] For our purposes, it will be convenient to work with a strengthened form of non-ordinarity. Recall that $f$ is $\ell$-non-ordinary if $a(\ell) \equiv 0 \pmod{\ell\mO_L}.$ In the case when $f$ is a normalized Hecke eigenform, $\ell$-non-ordinarity extends to congruences for many Fourier coefficients of $f$ modulo powers of $\ell.$

\begin{lemma}
\label{extendpowers}
Let $k\geq 4$ be an even integer and let $\ell$ be prime. Suppose $f(z) = \sum_{n = 0}^{\infty} a(n)q^n \in M_k \cap \mO_L[[q]]$ is a normalized Hecke eigenform. If $f$ is $\ell$-non-ordinary, then for all $r, n \ge 1$, \[ a(\ell^r n) \equiv 0 \pmod{\ell^r\mO_L}. \]
\end{lemma}

\begin{proof} 
We fix $n$ and prove the result by induction on $r$. Because $f$ is a normalized Hecke eigenform, we have 
\begin{equation*}
\label{normeigen}
    a(\ell n) = a(\ell)a(n) - \ell^{k - 1}a(n/\ell)
    \equiv 0 \pmod{\ell\mO_L}. \nonumber
\end{equation*}
This establishes the case $r = 1.$ 
For the inductive step, we apply the last equation with $n$ replaced by $\ell^rn$ to find \[ a(\ell^{r + 1}n) = a(\ell)a(\ell^rn) - \ell^{k - 1}a(\ell^{r - 1}n). \] 
By the induction hypothesis for $r$ and $r - 1$, we conclude that 
\[ a(\ell^{r + 1}n) \equiv 0 \pmod{\ell^{r + 1}\mO_L}.\qedhere \]
\end{proof}

It remains an open question \cite{story} whether a generic normalized Hecke eigenform is non-ordinary at an infinite number of primes.
Jin, Ma, and Ono \cite{non-ord} proved that if $S$ is a finite set of primes, there are infinitely many normalized Hecke eigenforms on $\SL_2(\Z)$ that are $\ell$-non-ordinary for each $\ell \in S.$ 
Their key result (Theorem 2.5 of \cite{non-ord}) identifies a sufficient condition to guarantee the $\ell$-non-ordinarity of all cusp forms with Fourier coefficients in $\mO_L.$ Combining Theorem 2.5 and Proposition 2.1 of \cite{non-ord}, we obtain the following result, which is crucial in the formulation of $\ell$ being good for $\aa.$

\begin{lemma}
\label{extraordinarycondition}
Let $k \ge 12$ be an even integer and let $f(z) = \sum_{n = 1}^{\infty} a(n)q^n \in S_k \cap \mO_L[[q]]$ be a normalized Hecke eigenform. Let $\ell$ be a prime such that $(\ell - 1) \mid (k - m)$ for some $m \in \{4, 6, 8, 10, 14\}.$ Then, $f$ is $\ell$-non-ordinary.
\end{lemma}

\begin{proof}
In Theorem 2.5 of \cite{non-ord}, we take integral $a \ge 0$ sufficiently large so that \[ k - 2 \le (m - 2)\ell^a. \] Clearly, $\ord_{\infty}(f) \ge 1 > -\ell^a$, so applying the theorem shows that \[ a(\ell^a) \equiv -\frac{2m}{B_m}a(0) \equiv 0 \pmod{\ell\mO_L}. \] By Proposition 2.1 of \cite{non-ord}, we conclude that $a(\ell) \equiv 0 \pmod{\ell\mO_L}$, i.e. $f$ is $\ell$-non-ordinary.
\end{proof}

Given an even weight $k \ge 12$, it is known that the space of cusp forms $S_k$ has a canonical basis of normalized Hecke eigenforms that is unique up to reordering. The coefficients of their Fourier expansions lie in the ring of integers of a number field $L_k.$ Thus, any cusp form $f \in S_k \cap \mO_{L_k}[[q]]$ can be written as an $L_k$-linear combination of normalized eigenforms.
If each normalized eigenform is $\ell$-non-ordinary, then $f$ will also be $\ell$-non-ordinary as long as the denominators of the coefficients in the linear combination are coprime to $\ell.$

To control these denominators, we introduce the \textit{weight $k$ Hecke determinant} for the space of cusp forms $S_k$, denoted by $\mD_k.$
Write $d_k := \dim S_k$, so that $S_k$ has a canonical basis of $d_k$ normalized eigenforms. 
We put the first $d_k$ coefficients of their Fourier expansions in each column of a $d_k \tm d_k$ matrix. We denote by $\mD_k \in \mO_{L_k}/\{\pm 1\}$ the determinant of this matrix. Because the basis of normalized eigenforms can be reordered, $\mD_k$ is only defined up to a sign.
Due to Cramer's rule, $\mD_k$ controls the denominators of the coefficients of a linear combination of the normalized eigenforms (see the proof of \Cref{goodpair}). For this reason, we are interested in primes $\ell$ for which the ideals $\mD_k\mO_{L_k}$ and $\ell\mO_{L_k}$ are coprime. This condition is equivalent to the more easily verifiable condition $\ell \nmid N_k(\mD_k)$, where $N_k$ is the norm over  the extension $L_k/\Q.$ When $\ell \nmid N_k(\mD_k)$, the $\ell$-non-ordinarity of normalized eigenforms extends through linear combinations.




\begin{lemma}
\label{goodpair}
Let $k \ge 12$ be an even integer and let $\ell$ be a prime such that $\ell \nmid N(\mD_k)$ and $(\ell - 1) \mid (k - m)$ for some $m \in \{4, 6, 8, 10, 14\}.$
Then for all $g = \sum_{n = 1}^{\infty} a_g(n)q^n \in S_k \cap \mO_{L_k}[[q]]$, and all positive integers $r$ and $n$, we have \[ a_g(\ell^r n) \equiv 0 \pmod{\ell^r\mO_{L_k}}. \] 
\end{lemma}

\begin{proof}
Let $f_1, \ldots, f_d$ be some ordering of the canonical basis of normalized Hecke eigenforms for $S_k.$ Write $L := L_k$ for the number field generated by their Fourier coefficients. By \Cref{extraordinarycondition}, each $f_i$ is $\ell$-non-ordinary. Write \[ f_i = \sum_{n = 1}^{\infty} a_i(n)q^n \in \mO_L[[q]]. \] There exist $\bb_i \in L_k$ such that \[ g = \sum_{i = 1}^d \bb_if_i. \] Therefore, we have the matrix equation \[ \bmat{a_1(1) & a_2(1) & \cdots & a_d(1) \\ a_1(2) & a_2(2) & \cdots & a_d(2) \\ \vdots & \vdots & \ddots & \vdots \\ a_1(d) & a_2(d) & \cdots & a_d(d)}\bmat{\bb_1 \\ \bb_2 \\ \vdots \\ \bb_d} = \bmat{a_g(1) \\ a_g(2) \\ \vdots \\ a_g(d)}. \] By Cramer's rule, we can write $\bb_i = \cc_i/\mD_k$ (where we fix some sign for $\mD_k$) for some $\cc_i \in \mO_L.$ Thus,
\begin{align}
    a_g(\ell^r n) &= \bb_1a_1(\ell^r n) + \cdots + \bb_da_d(\ell^r n) \nonumber \\
\label{agdiv}
    &= \frac{1}{\mD_k}\ps{\cc_1a_1(\ell^r n) + \cdots + \cc_da_d(\ell^r n)}.
\end{align}
Because each $f_i$ is $\ell$-non-ordinary, \Cref{extendpowers} implies that for all $i$, \[ a_i(\ell^r n) \equiv 0 \pmod{\ell^r\mO_L}. \] Because $\ell \nmid N(\mD_k)$, we know that $\mD_k\mO_L$ is coprime to $\ell\mO_L.$ Therefore, by \eqref{agdiv}, we conclude that $a_g(\ell^r n) \equiv 0 \pmod{\ell^r\mO_L}.$
\end{proof}

\subsection{Modular forms modulo $\ell$}
\label{modulo} 
To study modular forms modulo $\ell$, we make use of the Serre's theory of filtration. We collect the set of modular forms modulo $\ell$ of weight $k$ into the space \[ M_{k, \ell} := \{ f(z) \pmod{\ell} : f(z) \in M_k \cap \Z[[q]] \}. \] For an integer weight modular form $f$ on $\SL_2(\Z)$ with $\ell$-integral rational coefficients, we define the \textit{filtration of $f$ modulo $\ell$} by \[ \om_\ell(f) := \inf\{k : f(z) \pmod{\ell} \in M_{k, \ell}\}. \]


We conclude this section with two important lemmas used in the proof of \Cref{finitel}. The first lemma (Proposition 2.44 of \cite{web}) is a fundamental result that describes the effect of the Ramanujan Theta-operator on filtration. The $\Th$-operator is defined on power series in $q$ by 
\begin{align*}
    \Theta \ps{\sum_{n=0}^\infty a(n)q^n} := \sum_{n=0}^\infty na(n)q^n.
\end{align*}
In terms of differentials, we have $\Th = q\frac{d}{dq} = \frac{1}{2\pi i}\frac{d}{dz}.$

\begin{lemma}
\label{2.44}
If $\ell \ge 5$ is a prime and $f \in M_k \cap \Z[[q]]$, then $\Th(f) \pmod{\ell}$ is the reduction of a modular form modulo $\ell.$ Moreover, \[ \om_\ell(\Th f) = \om_\ell(f) + (\ell + 1) - s(\ell - 1) \] for some integer $s \ge 0$, with equality if and only if $\ell \nmid \om_\ell(f).$
\end{lemma}

The second lemma (Corollary of Proposition 2.56 of \cite{web}) controls the filtrations of powers of $\Delta$ acted on by the $\Th$-operator.

\begin{lemma}
\label{2.56}
If $\ell \ge 5$ is a prime and $\dd_\ell$ is a positive integer, then for any $m \ge 0$, \[ \om_\ell(\Th^m \Delta^{\dd_\ell}) \ge \om_\ell(\Delta^{\dd_\ell}) = 12\dd_\ell. \]
\end{lemma}

\section{Proofs of the Main Results}
\label{proofs}


\subsection{Proofs of \Cref{kiming} and \Cref{thm2}}

The key input for \Cref{kiming} is a detailed study of the $\Th$-operator on Serre filtrations, which Kiming and Olsson used to classify when a power of $\Delta$ is fixed by a power of $\Th$ \cite{K+O}. 

\begin{proof}[Proof of \Cref{kiming}]

Since $\ell \nmid b$, we may write $\aa = 24k + \ell u$ for some $k \ge 1$ and $u \in \Z_{(\ell)}$, where $\Z_{(\ell)}$ denotes the ring of $\ell$-integral rational numbers. By \Cref{movel}, we have
\begin{equation}
\label{deltaequiv}
    \sum_{n = 0}^{\infty} p_\aa(n)q^{n + k} = q^k\poc{q}^{24k + \ell u} = \Delta^k \poc{q}^{\ell u} \equiv \Delta^k\poc{q^\ell}^u \pmod{\ell}.
\end{equation}
Throughout the rest of this paper, we write $\Delta^k =: \sum_{n = 0}^{\infty} \tau_k(n)q^n.$ We extract terms of the form $q^{\ell n + c + k}$ on both sides of \eqref{deltaequiv} and use the assumption that $p_\aa(\ell n + c) \equiv 0 \pmod{\ell}$ for all $n$ to find \[ \tau_k(\ell n + c + k) \equiv 0 \pmod{\ell} \] for all $n.$ It follows from Fermat's little theorem that \[ \Th^{\ell - 1}\ps{q^{-(c + k)}\Delta^k} \equiv q^{-(c + k)}\Delta^k \pmod{\ell}. \] Applying Theorem 3 of \cite{K+O}, we conclude\footnote{As phrased in \cite{K+O}, this theorem can only be applied in our setting when $c + k \in \{0, 1, \ldots, \ell - 1\}.$ However, it is clear from the proof of the theorem that it applies for $c + k$ an arbitrary integer.} that $c + k \equiv 0 \pmod{\ell}.$ Consequently, \[ 24c + \aa \equiv 24(c + k) \equiv 0 \pmod{\ell}.\qedhere \]
\end{proof}

The proof of \Cref{thm2} makes use of Euler's and Jacobi's identities \eqref{euler} and \eqref{jacobi} for $\eta(24z)$ and $\eta(8z)^3$, respectively. We express the generating function for $p_\aa$ in terms of these formulae and compare coefficients to obtain congruences. This provides a small extension to Theorem 1.2 (1) and (2) from \cite{chan}.

\begin{proof}[Proof \Cref{thm2}]
Write $r := \ord_\ell(\aa-1)$ so that $\aa - 1 = \ell^r u$ for some $u \in \Z_{(\ell)}.$ It follows that \[ \sum_{n = 0}^\infty p_{\aa}(n)q^{24n + 1} = q\poc{q^{24}}^{1 + \ell^r u} = \eta(24z)\poc{q^{24}}^{\ell^r u}. \]
Substituting for $\eta(24z)$ with  Euler's identity \eqref{euler} and applying \Cref{movel}, we conclude that \[ \sum_{n=0}^\infty p_{\aa}(n)q^{24n + 1} \equiv \poc{q^{24\ell}}^{\ell^{r - 1}u} \sum_{n = 1}^{\infty} \leg{12}{n}q^{n^2} \pmod{\ell^r}. \]
The term $\poc{q^{24\ell}}^{\ell^{r - 1}u}$ is a power series in $q^{\ell}$, while the term $\sum_{n = 1}^{\infty} \leg{12}{n}q^{n^2}$ contains only monomials in $q$ whose exponent is a perfect square. Thus, if $24n + 1$ is not a square modulo $\ell$, the preceding equation implies that $p_{\aa}(n) \equiv 0 \pmod{\ell^r}.$ This proves part (i) of \Cref{thm2}. Part (ii) follows along the same lines by defining $r := \ord_\ell(\aa - 3)$ and using Jacobi's identity \eqref{jacobi}. 
\end{proof}

\subsection{Proof of \Cref{nonordthm}}
As explained in \Cref{intro}, \Cref{nonordthm} uses the $\ell$-non-ordinarity of modular forms to obtain $\ell^r$-balanced congruences for $p_\aa.$ To do so, we express the generating function for $p_\aa$ in terms of a power of Ramanujan's Delta-function. Under suitable hypotheses, this power of $\Delta$ is $\ell$-non-ordinary, from which we obtain $\ell^r$-balanced congruences for $p_\aa.$

\begin{proof}[Proof of \Cref{nonordthm}]

Write $r := \ord_\ell(24k -\alpha)$ so that $\aa - 24k = \ell^ru$ for some $u \in \Z_{(\ell)}.$
We prove by induction that for all $1 \le i \le r$, we have
\begin{equation}
\label{induct_thm2} 
    \sum_{n=0}^\infty p_{\alpha}(\ell^in - k)q^n \equiv \poc{q}^{\ell^{r - i}u}\sum_{n=0}^\infty \tau_k(\ell^in)q^n \pmod{\ell^r}. 
\end{equation} 
For the base case $i = 1$, we use \Cref{movel} to write
\begin{equation*}
    \sum_{n=0}^\infty p_{\alpha}(n)q^{n + k} = q^k\poc{q}^{24k + \ell^r u} = \Delta^k \poc{q}^{\ell^r u} \equiv \Delta^k \poc{q^\ell}^{\ell^{r - 1}u} \pmod{\ell^r}.
\end{equation*}
Extracting terms of the form $q^{\ell n}$ on both sides of the last equation and replacing $q^\ell$ by $q$, we obtain \[ \sum_{n=0}^\infty p_{\alpha}(\ell n - k)q^n \equiv \poc{q}^{\ell^{r - 1}u}\sum_{n=0}^\infty \tau_k(\ell n)q^n \pmod{\ell^r}. \] This establishes the base case.

For the inductive step, we  note that the conditions of \Cref{goodpair} are satisfied because $\ell$ is good for $\aa$ with parameter $k.$ Therefore, by \Cref{goodpair} we have \[ \tau_k(\ell^{i}n) \equiv 0 \pmod{\ell^{i}}. \] In addition, \Cref{movel} gives \[ \poc{q}^{\ell^{r - i}u} \equiv \poc{q^\ell}^{\ell^{r - (i + 1)}u} \pmod{\ell^{r - i}}. \] Thus, we can rewrite \eqref{induct_thm2} in the form \[ \sum_{n=0}^\infty p_{\alpha}(\ell^i n - k)q^n \equiv \poc{q^\ell}^{\ell^{r - (i+1)}u}\sum_{n=0}^\infty \tau_k(\ell^i n)q^n \pmod{\ell^r}. \] Extracting terms of the form $q^{\ell n}$ on both sides of the previous equation and replacing $q^\ell$ by $q$, we find \[ \sum_{n=0}^\infty p_{\alpha}(\ell^{i + 1}n - k)q^{n} \equiv \poc{q}^{\ell^{r - (i + 1)}u}\sum_{n=0}^\infty \tau_k(\ell^{i + 1}n)q^{n} \pmod{\ell^r}. \] This completes the induction.

To finish the proof, we take $i := v$ in \eqref{induct_thm2}. By \Cref{goodpair}, we know that $\tau_k(\ell^v n) \equiv 0 \pmod{\ell^v}$ for all $n.$ Thus, we conclude that \[ \sum_{n=0}^\infty p_{\alpha}(\ell^{v}n - k)q^n \equiv \poc{q}^{\ell^{r-v}u}\sum_{n=0}^\infty \tau_k(\ell^{v}n)q^n \equiv 0 \pmod{\ell^v}. \qedhere \] 
\end{proof}

\Cref{nonordthm} yields an $\ell^r$-balanced congruence for $p_\aa$ for all primes $\ell$ good for $\aa.$ In fact, for all but finitely many $\aa$, there are only finitely many primes $\ell$ good for $\aa.$ Furthermore, the associated parameter $k$ is strictly upper bounded by $\ell.$ 


\begin{proposition}\label{kltell}
If a prime $\ell$ is good for $\alpha$ with parameter $k$, then $k < \ell.$
\end{proposition}

\begin{proof}
Suppose for contradiction that $k \ge \ell.$ By \Cref{extraordinarycondition}, we know that all normalized Hecke eigenforms in $S_{12k}$ are $\ell$-non-ordinary. Therefore, each entry in the $\ell$th row of the matrix defining the weight $12k$ Hecke determinant $\mD_{12k}$ is divisible by $\ell.$ Consequently, $\ell \mid N_{12k}(\mD_{12k})$, in contradiction to \eqref{cond3}. This contradiction proves that $k < \ell.$
\end{proof}

We can use the preceding proposition to obtain an upper bound on a prime $\ell$ good for $\aa = a/b$ in terms of $a$ and $b.$ This proves that for all but finitely many $\aa$, \Cref{nonordthm} identifies only finitely many congruences for $p_\aa.$ 


\begin{proposition}
\label{finitegood}
Suppose that $\ell$ is good for $\alpha = a/b.$ If $\aa$ is not an even integer in the range $[-14, 26]$, then $\ell\leq \max\{\abs{26b - a}, \abs{14b + a}\}.$
\end{proposition}

\begin{proof}
Suppose that $\ell$ is good for $\aa$ with parameter $k.$ If $\ell \in \{2, 3\}$, the result is obvious. Thus, suppose that $\ell \ge 5.$ By \eqref{cond2}, we can write $12k - m = c(\ell - 1)$ for some $m \in \{4, 6, 8, 10, 14\}$ and $c \in \Z.$ Because $k \ge 1$, it is easy to see that $c \ge 1.$ By \Cref{kltell}, we have \[ 12(\ell - 1) \ge 12k = c(\ell - 1) + m. \] Thus, $c < 12.$ It follows that $c \le 11.$ We now use \eqref{cond1} to obtain
\begin{align*}
    0 \equiv b(24k - \aa) &= 24kb - a \\
    &= 2b(c(\ell - 1) + m) - a \\
    &\equiv 2b(m-c) - a \pmod{\ell}.
\end{align*}
Note that $-14 \le 2(m - c) \le 26.$ Because $\aa$ is not an even integer in the range$[-14, 26]$, we know that  $2b(m - c) - a \ne 0.$ It follows that $|2b(m - c) - a| \ge \ell.$ Testing the extremal values of $2(m - c)$, we obtain the desired result.
\end{proof}

\begin{remark*}
By checking the possible values of $2(m - c)$ and using \Cref{filtrationcor} to handle the case of $\aa < 0$, the preceding proposition can easily be extended to all $\aa \not \in \{0,2,4,6,8,10,14,26\}.$ This set of exceptional $\aa$ corresponds to the set of even lacunary powers of the eta-function classified by Serre.
\end{remark*}

\subsection{Proofs of \Cref{finitel} and \Cref{filtrationcor}}

For our next proof, we employ similar strategies to Ahlgren and Boylan \cite{ahlgren}, and Boylan \cite{boylan}, taking advantage of Serre filtrations of modular forms. We restate an $\ell$-balanced congruence for $p_\aa$ equivalently in terms of a modular form $f_\ell$ being fixed by $\Th^{\ell - 1}.$ We then study the sequence of weights $\om_\ell(\Th^if_\ell)$ to gain further information.


\begin{proof}[Proof of \Cref{finitel}]
Suppose for contradiction that $p_\aa$ admits an $\ell$-balanced congruence of the given form for some prime $\ell > |\aa| + 4.$  It follows that $\ell \ge 5$ and $\aa \not\equiv 0, 2, 4 \pmod{\ell}.$ Without loss of generality, we take $\dd_\ell$ to be a positive integer. By \Cref{kiming}, we know that $24\dd_ \ell \equiv \aa \pmod{\ell}.$ Let $24\dd_\ell = \alpha + \ell u$ for some $u \in \Z_{(\ell)}.$ It follows from \Cref{movel} that \[ \Delta^{\dd_\ell} = q^{\dd_\ell}\poc{q}^{\aa + \ell u} = \poc{q}^{\ell u}\sum_{n = 0}^{\infty} p_\aa(n - \dd_\ell)q^n \equiv \poc{q^\ell}^u\sum_{n=0}^\infty p_\aa(n - {\dd_\ell})q^n \pmod{\ell}. \]
Put $f_\ell := \Delta^{\dd_\ell}.$ Thus, the $\ell$-balanced congruence in the theorem statement is equivalent to \[ \Th^{\ell - 1} f_\ell \equiv f_\ell \pmod{\ell}. \] Let $c := \Psi_\ell(-12\dd_\ell).$ Because $-2c \equiv \aa \not\equiv 2, 4 \pmod{\ell}$, we have $c < \ell - 2.$ Repeatedly applying \Cref{2.44}, we observe that \[ \om_\ell(\Th^c f_\ell) = \om_\ell(f_\ell) + c(\ell + 1) \equiv 12\dd_\ell + c \equiv 0 \pmod{\ell}. \] By \Cref{2.44}, there exists an integer $s \ge 1$ such that 
\begin{align}
\label{c+1}
    \om_\ell(\Th^{c + 1}f_\ell) &= \om_\ell(\Th^c f_\ell)+ (\ell+1) - s(\ell - 1) \nonumber \\
    &= 12\dd_\ell + (c + 1)(\ell + 1) - s(\ell - 1).
\end{align}
Applying \Cref{2.56}, we deduce that \[ s \le \frac{(c + 1)(\ell + 1)}{\ell - 1} < c + 3. \]
Hence, $s \le c + 2.$ 
We now choose $j \ge 1$ minimal such that \[ \om_\ell(\Th^{c + j}f_\ell) \equiv 0 \pmod{\ell}. \] 
Such a $j$ exists and satisfies $j \le (\ell - 2) - c$ by the following argument. 
If $\om_\ell(\Th^{\ell - 2}f_\ell) \not\equiv 0 \pmod{\ell}$, then by \Cref{2.44}, we get \[ \om_\ell(\Th^{\ell - 2}f_\ell) = \om_\ell(\Th^{\ell - 1}f_\ell) - (\ell + 1) = \om_\ell(f_\ell) - (\ell + 1) < \om_\ell(f_\ell), \] in contradiction to \Cref{2.56}. 
Thus, it is imperative that $\om_\ell(\Th^{\ell - 2}f_\ell) \equiv 0 \pmod{\ell}.$ 
This proves that $j$ exists and $j \le (\ell - 2) - c.$  We now observe that \[ 0 \equiv \om_\ell(\Th^{c + j}f_\ell) = 12\dd_\ell + (c + j)(\ell + 1) - s(\ell - 1) \equiv j + s \pmod{\ell}. \] 
By the bounds on $j$ and $s$, we must have $(j, s) = (\ell - c - 2, c + 2).$ Rewriting \eqref{c+1}, we obtain
\begin{equation*}\label{weight}
   \om_\ell(\Th^{c + 1}f_\ell) = 12\dd_\ell + 2c - \ell + 3.  
\end{equation*}
By \Cref{2.56}, we deduce that 
\begin{equation}
\label{doublec}
    2c - \ell + 3 \ge 0.
\end{equation}
If $\aa < 0$ is even, then $2c = 2\Psi_\ell\ps{-\frac{\aa}{2}} = -\aa.$ Hence \[ 2c - \ell + 3 < -\aa - (|\aa| + 4) + 3 < 0. \] If $\aa > 3$ is odd, then $2c = 2\Psi_\ell\ps{-\frac{\aa}{2}} = \ell - \aa.$ Hence \[ 2c - \ell + 3 < (\ell - \aa) - \ell + 3 < 0. \] In either case, we derive a contradiction to \eqref{doublec}. This contradiction proves that $\ell \le |\aa| + 4.$
\end{proof}

\Cref{finitel} cannot be easily extended for non-integral $\aa$, because $c = \Psi_\ell\ps{-\frac{\aa}{2}}$ oscillates from $0$ to $\ell - 1$ as $\ell$ varies. However, we can use the proof of \Cref{finitel} for rational $\aa$ to constrain the set of primes $\ell$ for which $p_\alpha$ can admit an $\ell$-balanced congruence.

\begin{proof}[Proof of \Cref{filtrationcor}]
Write $\ell = 2bk + r$, where $r := \Psi_{2b}(\ell).$ Clearly $\gcd(r, 2b) = 1$, so we may define $s \in \Z$ such that \[ \Psi_{2b}\ps{\frac{a}{r}}r = 2bs + a. \] Put \[ c := \Psi_{2b}\ps{\frac{a}{r}}k + s. \] If $\aa$ is a negative even integer, then $\Psi_{2b}\ps{\frac{a}{r}} = 0$, and it is easy to check that $c = -\frac{\aa}{2} \ge 0.$ Otherwise $\Psi_{2b}\ps{\frac{a}{r}} \ge 1$, and hence \[ c \ge k + s \ge k - \frac{a}{2b} = \frac{\ell - (r + a)}{2b} \ge 0. \] Furthermore, we check that \[ 2bc = \Psi_{2b}\ps{\frac{a}{r}}(2bk) + 2bs \equiv \Psi_{2b}\ps{\frac{a}{r}}(-r) + 2bs = -a \pmod{\ell}. \] Thus, $c \equiv -\frac{\aa}{2} \pmod{\ell}.$ Suppose for contradiction that $\Psi_{2b}\ps{\frac{a}{r}} < b.$ It follows that \[ 2s - r = \frac{\Psi_{2b}\ps{\frac{a}{r}}r - a}{b} - r < -\frac{a}{b}. \] Therefore, we estimate 
\begin{align}
\label{2cestimate}
    2c - \ell + 3 &= k\ps{2\Psi_{2b}\ps{\frac{a}{r}} - 2b} + (2s - r) + 3 \nonumber \\
    &< -2k - \frac{a}{b} + 3 \nonumber \\
    &= \frac{(3b + r - a) - \ell}{b} < 0.
\end{align}
In particular, $c < \ell$ so $c = \Psi_\ell\ps{-\frac{\aa}{2}}.$ Because $\ell \ge |a| + 5b > |\aa| + 4$ and $\aa \not\in \{0, 2, 4\}$, the proof of \Cref{finitel} ensures that \eqref{doublec} still holds. However, \eqref{2cestimate} contradicts \eqref{doublec}. This contradiction proves that if $p_\aa$ admits an $\ell$-balanced congruence for $\ell \ge |a| + 5b$, then $\Psi_{2b}(\frac{a}{r}) \ge b$. 
\end{proof}


\section*{Acknowledgements}

The authors would like to thank Ken Ono, Larry Rolen, and Ian Wagner for their guidance on this project. This research was generously supported by the Asa Griggs Candler Fund, National Security Agency Grant H98230-19-1-0013,
National Science Foundation Grants 1557960 and 1849959, the Spirit of Ramanujan Global STEM Talent Search, and Princeton University.


\medskip
\begin{footnotesize}
  E. Bevilacqua, \textsc{Department of Mathematics, Penn State University,
    University Park, PA 16802}\par\nopagebreak
  \textit{E-mail address:} \texttt{eib5092@psu.edu}

  \medskip

  K. Chandran, \textsc{Department of Mathematics,
    Princeton University, Princeton, NJ 08544}\par\nopagebreak
  \textit{E-mail address:} \texttt{kapilchandran@princeton.edu}

  \medskip

  Y. Choi, \textsc{Phillips Exeter Academy,
    Exeter, NH 03833}\par\nopagebreak
  \textit{E-mail address:} \texttt{ychoi@exeter.edu}
\end{footnotesize}
  
\end{document}